\documentclass{svproc}
\usepackage{url}
\usepackage{upgreek}
\usepackage{amsmath}
\usepackage{amssymb}
\usepackage{xfrac}
\usepackage{subcaption}
\usepackage{pgfornament}
\usepackage{tikz}

\usepackage{csquotes}
\usepackage{xcolor}
\definecolor{MapleRed}{rgb}{1,0,0}
\definecolor{MapleBlue}{rgb}{0,0,1}
\definecolor{MaplePink}{rgb}{1,0,1} 

\def\MapleInput#1{\noindent{{\small $>$ {\tt \color{MapleRed}{#1} }}}}
\def\MapleOutput#1{\textsl{\begin{center} \color{MapleBlue}{#1} \end{center}}}

\begin{document}
\mainmatter              
\title{A Random Walk through Experimental Maths}
\titlerunning{Experimental Maths}  
%
\author{Eunice Y.~S.~Chan\inst{1} \and Robert M.~Corless\inst{1}}
\authorrunning{Chan and Corless} 
%
\tocauthor{Eunice Y.~S.~Chan, Robert M.~Corless}
\institute{ORCCA and SMSS and The Rotman Institute of Philosophy\\ Western University, London, N6A 5B7, Canada \\
\email{echan295@uwo.ca}, \email{rcorless@uwo.ca}}

\maketitle              

\begin{abstract}
We describe our adventures in creating a new first-year course in Experimental Mathematics that uses active learning. We used a state-of-the-art facility, called The Western Active Learning Space, and got the students to `drive the spaceship' (at least a little bit). This paper describes some of our techniques for pedagogy, some of the vignettes of experimental mathematics that we used, and some of the outcomes. EYSC was a student in the simultaneously-taught senior sister course ``Open Problems in Experimental Mathematics" the first time it was taught and an unofficial co-instructor the second time. Jon Borwein attended the Project Presentation Day (the second time) and gave thoughtful feedback to each student. This paper is dedicated to his memory.
\end{abstract}
\section{Introduction}\label{s:1}
It has been known at least since the nineties that active learning is the most effective method for teaching science \cite{handelsman2004scientific}. It is also clear to many people that the corpus of mathematical knowledge has changed dramatically since the introduction of computers. It is less widely acknowledged that this ought to change not only \textsl{how} we teach mathematics, but \textsl{what} mathematics we should teach in the first place.

Of course, there \textsl{have} been vocal advocates of just this, including for instance Gilbert Strang \cite{strang2001too} and, very notably, Jon Borwein and co-workers. One of the present authors has held this view for decades, and put it into practice on several occasions, see e.g. \cite{corless2004computer}, \cite{corless1995first}, and \cite{corless1997scientific}.

A major problem, pointed out by Hamming in his iconoclastic Calculus textbook~\cite{hamming2012methods}, is that making only minor changes to our practice can only make things worse, if we're at a local optimum. To really do better, we have to make very large changes, basically all at once.

The difficulties with \textsl{that} are probably obvious, but we'll list some of them anyway: first, resources (it takes time and effort to design from scratch, not to mention to learn enough to use best pedagogical practices); second, inertia (the standard calculus, linear algebra, higher math sequence has enormous installed infrastructure including expectations); and finally outright hostility to change, amounting to denial.

In this paper, we'll talk about how we took advantage of an unusual confluence of opportunities to make a serious attempt, serious enough that we feel that this description will be useful for the next attempt. 

\subsection{The confluence of opportunities}
In May 2014, David Jeffrey (then Chair of Applied Mathematics at Western) was lamenting the invisibility of our program to students in their first year of study. A decade previously, you see, the first year calculus and linear algebra offerings of the separate Mathematics Department and of ours had been merged, so students only saw ``Calculus" or ``Linear Algebra." The Applied Math undergraduate program slowly declined, afterwards. In response to the Chair's lament, RMC simply wrote ``Introduction to Experimental Mathematics"---just that---on the whiteboard. The Chair and subsequently surveyed students were very enthusiastic, and so, with a verbal promise of a two-year trial both from the Chair and from the Dean of Science, RMC designed the course over June, July and August, making (with undergraduate student help) the recruiting video \url{https://vimeo.com/itrc/review/99140780/68995faea3} and taught the newly-christened AM 1999 September--December, of that same year, 2014.

With such short notice, enrolment was going to be an issue---Western does not give teaching credit for undergraduate courses with 9 or fewer students, or graduate courses with 4 or fewer students. Therefore, RMC invented and designed the sister course AM 9619 intended for graduate students and senior undergraduate students (such as EYSC, then---she's now a PhD student of RMC's at this time of writing). Again a verbal promise of a two-year trial was given by the Chair and Dean. These were to be taught simultaneously, and credit given for one course taught.
\begin{displayquote}
	``[on seeing the video] This is a wonderful course! It should be profiled!" --- Charmaine Dean, Dean of Science
\end{displayquote}

The major resources available for this were first RMC's time---he used a (northern) summer, normally used for research, to design the course, make the video, socialize it with students, design the grad course, attend the ICERM conference on Experimental Mathematics organized by Jon Borwein and David Bailey to seek advice there---Neil Calkin was especially helpful---explain the course to academic counselors, make a poster for Fall Preview Day, etc. Another task was to formally write up the course design and get it approved by Senate. Later, the Course Learning Outcomes had to be drafted for the undergraduate program description and harmonized with the Program Learning Outcomes. Second, RMC had some discretionary money he had available to pay the first TA, Steven Thornton, and an assistant, Torin Viger, to construct course materials (and to travel to the ICERM conference). The third resource was the \$400,000 Western Active Learning Space, which had just been built (RMC gave the very first class in it). We'll talk more about this space and how we used it later. Details of the space can be found at \url{http://www.uwo.ca/wals/}. The fourth resource was RMC's record in Experimental Mathematics and in the use of technology for teaching mathematics, which helped sell the course (both to the students and to the higher administration).

The members of the Undergraduate Society of Applied Mathematics (USAM) constituted a crucial fifth resource. They were consulted frequently in the ``summer of design'' and many of them helped to make the recruitment video linked above. Their enthusiastic support for the course was very heartening.

\subsection{The Western Active Learning Space}
The advanced features of WALS that we utilized included, first, the ``Air Media'' by which students or the TA could share their computer screens with their 7-member pods or with the whole class. This feature is not so expensive and one could imagine classes less richly furnished that still had it. WALS had multiple screens (one per pod) that could be used independently, and this was more expensive (but still useful). We also used the document camera frequently. However, we did not utilize the classroom's sound system or SMART Board pens.

We did use the low-tech whiteboards---there were three sizes: main boards, portable boards in A-frames and small ``huddle boards'' that could be placed at will on the main boards.

WALS came with eight recent vintage laptops each equipped with Maple. We used these each class; they were intended to be one per pod, plus the instructor's station, but we didn't always use them that way. The senior students had Maple on their personal laptops but the undergraduates did not. Therefore, priority of the classroom laptops was given to the undergraduate students.

\subsection{Active learning techniques}
Active learning techniques run from the obvious (get students to choose their own examples, and share) through the eccentric (interrupt students while programming similar but different programs and have them trade computers and problems) to the flaky (get them to do an interpretive dance or improvisational skit about their question). We tried to avoid the extremely flaky, but we did mention them, so that these introverted science students knew that this was within the realm of possibility.

The simplest activity was typing Maple programs that were handwritten on a whiteboard into a computer: this was simple but helpful because students learned the importance of precision, and had \textsl{immediate} help from their fellow students and from the TA.

Next in complexity was interactive programming exercises (integrated into the problems). Mathematicians tend to under-value the difficulty of learning syntax and semantics simultaneously.

We describe our one foray into eccentricity. The paper Strange Series and High Precision Fraud by Borwein and Borwein \cite{borwein1992strange} has six similar sums. We had six teams program each sum, at a stage in their learning where this was difficult (five weeks into the course). After letting the teams work for twenty minutes, we forced one member of each team to join a new team; each team had to explain their program (none were working at this stage) to the new member. This exercise was most instructive. The lessons learned included:
\begin{itemize}
	\item people approach similar problems very differently
	\item explaining what you are doing is as hard as doing it (maybe harder)
	\item basic software engineering (good variable names, clear structure, economy of thought) is important
	\item designing on paper first might be a good idea (nobody believed this, really, even after)
	\item social skills matter (including listening skills).
\end{itemize}
Perhaps the most important ``active learning" technique used, and the hardest to describe accurately, was a shift in perspective: the instructor attempted to listen carefully to students' wishes, and alter the activity, discussion, or topic depending on what their questions were. Waiting long enough for students to actually ask questions was often uncomfortable for everyone. The students were, at least initially, reluctant to use their freedom.

\subsection{Active learning in mathematics}
It is widely accepted that ``true'' or ``deep'' learning only happens when students are actively engaged. The old saying goes, ``I hear, and I forget. I see, and I remember. I do, and I understand.'' There is a lot of argument as to how best to make the students active, engaged, and make them do the work. There is a lot of discussion of the value of ``discovery'' versus the apparent efficiency of simple delivery of a lecture, or reading. We won't settle those arguments here because in our course, a mix was employed: some actual lectures (most often short, usually explaining how to program something in Maple), some individual activities, some partner work, some games, some peer assessment, a structured project, and giving them choices. RMC found it very hard (sometimes) to let the students make mistakes. There was an improvisational element: every class had a plan, and a goal (a ``planned learning outcome'') even if only a modest one, but at any moment tangential discussions could be seized as opportunities.

Another precious resource for this course was that it was not a prerequisite for anything---it had no required specific topics to cover---and thus the learning outcomes could be quite general, such as achieving the precision of thought necessary to program computers. Therefore, the class could afford to sail off on tangents, chosen by members of the class.

The students found this freedom (``lack of structure'') frightening at first, but exhilarating at the same time. They learned to trust the instructor, and the instructor learned to trust them.

This is not to say that there was \textsl{no} structure: there certainly was. The main theme of the course was that of discrete dynamical systems, normally taught long after calculus and analysis. But one of the major achievements of computers for mathematics research (and education) is in making deep questions immediately accessible to students, questions such as ``which initial guesses for Newton's method converge to which roots?'' This question immediately gives access to fractals, and their compelling images.

It is less appreciated that computers, via symbolic computation programs such as Maple or Mathematica, give students the same (or greater!) computational power to perform mathematical experiments that the giants of mathematical history had; Lagrange, Stirling, or even Euler, Gauss, and Newton~\cite{borwein2009computer}. The students can explore classical topics such as continued fractions experimentally and make deep discoveries themselves.

We ended the course with the Chaos Game Representation of DNA sequences. This has recently been used to generate an objective ``Map of Life''~\cite{kari2015mapping}. We point out that at the time the course was first taught, this was new (not yet published) frontier research. The students understood it. We got them from entering university to front-line research results. The rest of this paper gives some examples of what we did, to do so.

One final introductory point: what did the students say they needed from our program, that we weren't delivering? The number one request (in our survey results, prior to designing the course) was ``more programming.'' The number two consideration was ``more say in the syllabus,'' that is, student-centred curriculum design. We aimed to address both needs with this course.

\section{Choices of Topics}
As previously stated, the theme of the course was ``discrete dynamical systems.'' This was chosen because it benefits massively from computer support,  it's outside the normal curriculum (relieving pressure to teach to a specific goal), and it's accessible, beautiful, interesting, has applications, and because it has deep connections to classical topics in number theory (itself outside our standard curriculum). Before starting to teach, RMC drew up a list of potential vignettes. The top of the list was the classical theory of simple continued fractions. This topic is extremely accessible (we will give the first vignette in Section~\ref{sec:cf} below) and at the same time interesting, deep, and demanding of introductory programming skills (iteration, conditionals, and induction for proof of correctness). The next was the solution of nonlinear equations by Newton's method and its generalizations. Other topics included the Thue-Morse sequence, the Online Encyclopedia of integer sequences, computation of $\pi$, the game of life, the numerical solution of differential equations by Euler's method, visualization of complex functions by phase plots, and chaos game representation of DNA sequences. We sample from these below.

The graduate course AM 9619 had an extra lecture per week (the senior students attended the AM 1999 activities together with the first year students, but they got some advanced material every week that the AM 1999 students did not). We will give one sample below.

\subsection{Continued Fractions and Rational approximations of $\sqrt{2}$}\label{sec:cf}
How does one give control over the pedagogy to the students? RMC could only think to do it gradually, by doing an example and then asking students to \textsl{choose} their own examples which we would (first) do together and then (for other choices) they would do themselves. The students were presented with the sequence
\begin{equation}
	1 \>, \dfrac{3}{2} \>, \dfrac{17}{12} \>, \dfrac{577}{408} \>, \dfrac{665857}{470832} \>, \>\cdots \>,
\end{equation}
which at that point was plucked from thin air. Each term $x_n$ is generated from its predecessor\footnote{Here, $x_0 = 1$, $x_1 = \sfrac{3}{2}$, so on and so forth.} by the rule $x_n = \sfrac{\left(x_{n-1} + \sfrac{2}{x_{n-1}}\right)}{2}$. What means the same thing, if we label the numerators and denominators by $x_n = \sfrac{p_n}{q_n}$, in other words
\begin{align}
	p_n &= p_{n-1}^2 + 2q_{n-1}^2 \nonumber \\
	q_n &= 2p_{n-1}q_{n-1} \>.
\end{align}
\begin{displayquote}
``At first glance, nothing seems simpler or less significant than writing a number, say $\sfrac{9}{7}$, in the form
\begin{equation}
	\dfrac{9}{7} = 1 + \dfrac{2}{7} = 1 + \cfrac{1}{\sfrac{7}{2}} = 1 + \dfrac{1}{3 + \sfrac{1}{2}} = 1 + \cfrac{1}{3 + \cfrac{1}{1 + \sfrac{1}{1}}} \>.
\end{equation}
It turns out, however, that fractions of this form, called \textsl{continued fractions}, provide must insight...'' --- from p.~3 of C.~D.~Olds, ``Continued Fractions"~\cite{olds1963continued}.
\end{displayquote}
Carl Douglas Olds won the 1973 Chauvenet Prize, the highest award for mathematical exposition, for his paper ``The Simple Continued Fraction for $e$." The book cited above is likewise a model of lucidity, and reads very well today. A new book \cite{borwein2014neverending} is similarly valuable to students.

What follows is a simulation of a whiteboard discussion.
\begin{center}
	\pgfornament[scale=0.5]{88}
\end{center}
What's happening [in Olds' example]? Let's do the same thing with each $x_n$. First, we take out the integer part. For our first two numbers, nothing much happens:
\begin{align}
	x_0 &= 1\\
	x_1 &= \dfrac{3}{2} = 1 + \dfrac{1}{2} = 1 + \cfrac{1}{1 + \sfrac{1}{1}} \>,
\end{align}
but this last isn't much use.

The next number is more interesting:
\begin{align}
	x_2 &=\dfrac{17}{12} = \dfrac{12 + 5}{12} = 1 + \dfrac{5}{12} \nonumber \\
	&= 1 + \cfrac{1}{\sfrac{12}{5}} = 1 + \cfrac{1}{2 + \sfrac{2}{5}} = 1 + \cfrac{1}{2 + \cfrac{1}{\sfrac{5}{2}}} \nonumber \\
	&= 1 + \cfrac{1}{2 + \cfrac{1}{2 + \sfrac{1}{2}}} \>.
\end{align}
The crucial step in this process is writing the fractional part that we get, after taking out the integer part, as a reciprocal of another fraction, i.e.:
\begin{equation}
	\dfrac{5}{12} = \cfrac{1}{\sfrac{12}{5}} \>.
\end{equation}
	
Now a longer example:
\begin{alignat}{5}
	x_3 = \dfrac{577}{408} &= \dfrac{408 + 169}{408} & &= 1 + \dfrac{169}{408} \qquad = 1 + \dfrac{1}{\sfrac{408}{169}} \nonumber \\
	&= 1 + \cfrac{1}{2 + \sfrac{70}{169}} & &= 1 + \cfrac{1}{2 + \cfrac{1}{\sfrac{169}{70}}} \nonumber \\
	&= 1 + \cfrac{1}{2+\cfrac{1}{2 + \sfrac{29}{70}}} & &= 1 + \cfrac{1}{2 + \cfrac{1}{2 + \cfrac{1}{\sfrac{70}{29}}}} \nonumber \\
	&= 1 + \cfrac{1}{2 + \cfrac{1}{2 + \cfrac{1}{2 + \sfrac{12}{19}}}} & &= 1 + \cfrac{1}{2 + \cfrac{1}{2 + \cfrac{1}{2 + \cfrac{1}{\sfrac{29}{12}}}}} \nonumber \\
	&= 1 + \cfrac{1}{2 + \cfrac{1}{2 + \cfrac{1}{2 + \cfrac{1}{2 + \sfrac{5}{12}}}}} & &= 1 + \cfrac{1}{2 + \cfrac{1}{2 + \cfrac{1}{2 + \cfrac{1}{2 + \cfrac{1}{\sfrac{12}{5}}}}}} \nonumber \\
	&= 1 + \cfrac{1}{2 + \cfrac{1}{2 + \cfrac{1}{2 + \cfrac{1}{2 + \cfrac{1}{2 + \sfrac{2}{5}}}}}} & &= 1 + \cfrac{1}{2 + \cfrac{1}{2 + \cfrac{1}{2 + \cfrac{1}{2 + \cfrac{1}{2 + \cfrac{1}{2 + \sfrac{1}{2}}}}}}} \nonumber \\
	&= 1 + [2 \>, 2 \>, 2\>, 2\>, 2 \>, 2 \>, 2] \quad \text{for short.}
\end{alignat}
At this point, you may feel like sticking out your tongue and giving us a raspberry for such obvious cheating. Think of it like ``television wrestling" and give the entertainment a chance!

When you think about it, it \textsl{is} a bit mysterious that the simple rule
\begin{equation}
	x_n = \dfrac{\left(x_{n-1} + \sfrac{2}{x_{n-1}}\right)}{2}
\end{equation}
can generate the continued fractions
\begin{equation}
	1 \>, 1 + [2] \>, 1 + [2\>, 2\>, 2]\>, \text{and } 1 + [2\>, 2\>, 2\>, 2\>, 2\>, 2\>, 2] \>.
\end{equation}
The next one,
\begin{equation}
	x_4 = \dfrac{665857}{470832} = 1 + [2\>, 2\>, 2\>, 2\>, 2\>, 2\>, 2\>, 2\>, 2\>, 2\>, 2\>, 2\>, 2\>, 2\>, 2]
\end{equation}
has fifteen 2's in it. (By the way; don't worry, we'll check that by computer later.) That's one, three, seven, and fifteen twos. What's next? We'll leave that for now and go back to the first question, about $x_n^2 = \sfrac{p_n^2}{q_n^2}$.

The squares of our sequence are
\begin{multline}
1 \>, \dfrac{9}{4} = 2\dfrac{1}{4} \>, \left(\dfrac{17}{12}\right)^2 = \dfrac{289}{144} = \dfrac{288 + 1}{144} = 2 + \dfrac{1}{144} = 2 + \dfrac{1}{12^2} \>,\\
\left(\dfrac{577}{408}\right)^2 = \dfrac{332929}{166464} = \dfrac{332928 + 1}{166464} = 2 + \dfrac{1}{166464} = 2 + \dfrac{1}{408^2}
\end{multline}
and at this point, we might be prepared to bet that 
\begin{equation}
	x_4^2 = \left(\dfrac{665857}{470832}\right)^2 = 2 + \dfrac{1}{470832^2} \doteq 2 + 4.5\times10^{-12} \>.
\end{equation}
Checking using RMC's phone (a Nexus 5), we see that this is, in fact, true. But what does it mean?

One thing it means is that our sequence can be written as
\begin{equation}
	\sqrt{2 - \dfrac{1}{1^2}} \>, \sqrt{2 + \dfrac{1}{2^2}} \>, \sqrt{2 + \dfrac{1}{12^2}} \>, \sqrt{2 + \dfrac{1}{408^2}} \>, \sqrt{2 + \dfrac{1}{470832^2}} \doteq \sqrt{2 + 4.5\times10^{-12}} \>,
\end{equation}
that is a sequence of square roots of numbers that rapidly approach 2. The denominator of $x_5$ is
\begin{equation}
	q_5 = 2p_4q_4 = 2\cdot 470832 \cdot 665857 \doteq 6.5\times10^{11} \>;
\end{equation}
the next
\begin{equation}
	\left(\dfrac{p_5}{q_5}\right)^2 = 2 + \dfrac{1}{q_5^2} \doteq 2 + 2\times10^{-24} \>,
\end{equation}
about as close to 2 as one molecule in a mole\footnote{Avogadro's number is $6\cdot 10^{23}$, about.}.
\begin{center}
	\pgfornament[scale=0.5]{88}
\end{center}
Here ends our simulated whiteboard discussion. Some more questions present themselves. Does this continue? Is $x_5 = 1 + [2 \>, 2 \>, \ldots, 2]$ with thirty-one 2's in the continued fraction? Does $x_6$ have sixty-three 2's in it? Is $x_n^2 = 2 + \sfrac{1}{q_n^2}$ always? Does this mean that $x_n \doteq \sqrt{2}$?

Here's another question. What is
\begin{equation}
	1 + \cfrac{1}{2 + \cfrac{1}{2 + \cfrac{1}{2 + \cfrac{1}{2 + \substack{\ \ \\ \\ \ddots}}}}}
\end{equation}
where the 2's continue forever? Does this make sense? At this point, many students were surprised at the perfect predictability, and repeating nature, of this continued fraction, because it is indeed true that with quite natural definitions, this infinite continued fraction can only be $\sqrt{2}$.

But ``everybody knows'' that the decimal expansion for $\sqrt{2}$ does not repeat, because $\sqrt{2}$ is irrational! Why is this different? Is it something special about $\sqrt{2}$? (Of course a continued fraction is not a decimal expansion.)

The students were then asked to choose their own examples. Predictably, they chose only other square roots ($\sqrt{3}$, $\sqrt{5}$, etc.~) until prodded.

After enough hand calculation, Maple was introduced as a calculator. Some hand calculation is necessary, at first, because the students need to \textsl{feel} that they are connected to the mathematics; they need to own it. Once they feel ownership (and once the computations get tedious) then computer assistance is welcomed.

Over the next three classes, the pace gradually increased until the Gauss map $G: x\mapsto \text{frac}(\sfrac{1}{x})$ mapping $(0, 1)$ to $[0, 1)$ seemed natural; patterns were identified by the students, such as termination implies rationality, ultimate periodicity implies that the number is a quadratic irrational. This is Lagrange's classical result, and not easy to prove. No proofs are given in this section, only computation; Galois' result that purely periodic continued fractions arise from \textsl{reduced} quadratic irrationals is mentioned but not emphasized---we concentrate on the interval $[0, 1)$.

By this time the students are engaged, talking to each other, choosing different numbers. The number $e$ causes amazement with its pattern: $e = 2 + [1, 2, 1, 1, 4, 1, 1, 6, 1, 1, 8, \ldots]$. The fact that the continued fraction for $\pi$ is not known and has no apparent pattern, causes astonishment: in the first two weeks we have arrived at a deep open problem.

By now the students are hooked. They've learned some Maple; they've learned some mathematics. They've learned about fixed points (equal points of period 1, including the golden ratio $\phi = 1 + [1, 1, 1, \ldots]$) and about periodic points. the journey has begun.

\subsection{Newton's method, $\sqrt{2}$, and solving nonlinear equations}
In the first vignette, we met with the sequence
\begin{equation}
	1 \>, \dfrac{3}{2} \>, \dfrac{17}{12} \>, \dfrac{577}{408} \>, \dfrac{665857}{470832} \>, \ldots
\end{equation}
which was generated by $x_{n+1} = \sfrac{\left(x_{n} + \sfrac{2}{x_n}\right)}{2}$; in words, the average of the number and two divided by the number. This vignette explores where that sequence came from, and its relationship to $\sqrt{2}$. We approached this algebraically, as Newton did. Consider the equation
\begin{equation}
	x^2 - 2 = 0 \>.
\end{equation}
Clearly the solutions to this equation are $x = \sqrt{2}$ and $x = -\sqrt{2}$. Let us \textsl{shift the origin} by putting $x = 1 + s$; so $s = 0$ corresponds to $x = 1$. Then
\begin{equation}
	\left(1 + s\right)^2 - 2 = 1 + 2s + s^2 - 2 = -1 + 2s + s^2 = 0 \>.
\end{equation}
we now make the surprising assumption that $s$ is so small that we may ignore $s^2$ in comparison to $2s$. If it turned out that $s = 10^{-6}$, then $s^2 = 10^{-12}$, very much smaller than $2s = 2\cdot10^{-6}$; so there are small numbers $s$ for which this is true; but we don't know that this is true, here. We just hope.

Then if $s^2$ can be ignored, our equation becomes
\begin{equation}
	-1 + 2s = 0
\end{equation}
or $s = \sfrac{1}{2}$. This means $x = 1 + s = 1 + \sfrac{1}{2} = \sfrac{3}{2}$.

We now repeat the process: shift the origin to $\sfrac{3}{2}$, not $1$: put now
\begin{equation}
	x = \sfrac{3}{2} + s
\end{equation}
so
\begin{equation}
	\left(\sfrac{3}{2} + s\right)^2 = \sfrac{9}{4} + 3s + s^2 - 2 = 0 \>.
\end{equation}
This gives $3s + s^2 + \sfrac{1}{4} = 0$ and again we ignore $s^2$ and hope it's smaller than $3s$. This gives
\begin{equation}
	3s + \sfrac{1}{4} = 0
\end{equation}
or $s = -\sfrac{1}{12}$. This means $x = \sfrac{3}{2} - \sfrac{1}{12}$ or $x = \sfrac{17}{12}$. Now we see the process. Again, shift the origin: $x = \sfrac{17}{12} + s$. Now
\begin{equation}
	\left(\dfrac{17}{12} + s\right)^2 = \dfrac{289}{144} + \dfrac{17}{6}s + s^2 - 2 = 0 \>.
\end{equation}
Ignoring $s^2$,
\begin{equation}
	\dfrac{17}{6}s + \dfrac{1}{144} = 0
\end{equation}
or
\begin{equation}
	s = \dfrac{-6}{17\cdot144} = \dfrac{-1}{17\cdot24} = \dfrac{-1}{408} \>.
\end{equation}
Thus,
\begin{equation}
	x = \dfrac{17}{12} - \dfrac{1}{408} = \dfrac{577}{408} \>.
\end{equation}
As we saw in the previous vignette, there are the exact square roots of numbers ever more close to 2. For instance,
\begin{equation}
	\dfrac{577}{408} = \sqrt{2 + \dfrac{1}{408^2}} \>.
\end{equation}
It was Euler who took Newton's ``shift the origin" strategy and made a general method---that we call Newton's method---out of it. In modern notation, Euler considered solving $f(x) = 0$ for differentiable function $f(x)$, and used the tangent line approximation near an initial guess $x_0$: if $x = x_0 + s$ then, using $f'(x_0)$ to denote the slope at $x_0$, $0 = f(x) = f(x_0 + s) \doteq f(x_0) + f'(x_0)s$ ignoring terms of order $s^2$ or higher. Then
\begin{equation}
	s = -\dfrac{f(x_0)}{f'(x_0)}
\end{equation}
so
\begin{equation}
	x \doteq x_0 + s = x - \dfrac{f(x_0)}{f'(x_0)} \>.
\end{equation}
The fundamental idea of Newton's method is that, if it worked once, we can do it again: pass the parcel! Put
\begin{align}
	x_1 &= x_0 - \dfrac{f(x_0)}{f'(x_0)} \\
	x_2 &= x_1 - \dfrac{f(x_1)}{f'(x_1)} \\
	x_3 &= x_2 - \dfrac{f(x_2)}{f'(x_2)}
\end{align}
and keep going, until $f(x_k)$ is so small that you're happy.

Notice that each $x_k$ solves
\begin{equation}
	f(x) - f(x_k) = 0
\end{equation}
not $f(x) = 0$. But if $f(x_k)$ is really small, you've solved ``almost as good" an equation, like finding $\sqrt{2 + \sfrac{1}{408^2}}$ instead of $\sqrt{2}$. So where did $\sfrac{\left(x_n + \sfrac{2}{s_n}\right)}{2}$ come from?
\begin{equation}
	x_{n+1} = x_n - \dfrac{f(x_n)}{f'(x_n)} = x_n - \dfrac{\left(x_n^2 - 2\right)}{2x_n}
\end{equation}
because if $f(x) = x^2 - 2$, $f'(x) = 2x - 0 = 2x$. Therefore,
\begin{alignat}{4}
	x_{n+1} &= x_n - \dfrac{\left(x_n^2 - 2\right)}{2x_n} & &= \dfrac{2x_n^2 - x_n^2 + 2}{2x_n} \nonumber \\
	&= \dfrac{x_n^2 + 2}{2x_n} & &= \dfrac{1}{2}\left(x_n + \dfrac{2}{x_n}\right)
\end{alignat}
as claimed\footnote{For more complicated functions one \textsl{shouldn't} simplify for numerical stability reasons. But for $x^2 - a$, it's okay.}.

Executing this process in decimals, using a calculator (our handy HP48G+ again), with $x_0 = 1$, we get
\begin{align}
	x_0 &= 1 \nonumber \\
	x_1 &= \underline{1}.5 \nonumber \\
	x_2 &= \underline{1.4}1666\ldots \nonumber \\
	x_3 &= \underline{1.41421}568628 \nonumber \\
	x_4 &= \underline{1.41421356238} \nonumber \\
	x_5 &= x_4 \text{\ to all 11 places in the calculator}
\end{align}
Now $\sqrt{2} = 1.41421356237$ on this calculator. We see (approximately) 1, 2, 5 then 10 correct digits. The convergence behaviour is clearer in the continued fraction representation:
\begin{multline}
	1 \>, 1 + \left[2\right] \>, 1 + \left[2 \>, 2 \>, 2\right] \>, 1 + \left[2 \>, 2\>, 2\>, 2\>, 2\>, 2\>, 2 \right]\>, \\
	1 + \left[2 \>, 2\>, 2\>, 2\>, 2\>, 2\>, 2\>, 2\>, 2\>, 2\>, 2\>, 2\>, 2\>, 2\>, 2\right]
\end{multline}
with 0, 1, 3, 7, 15 twos in the fraction part: each time doubling the previous plus 1, giving $2^0 - 1$, $2^1 - 1$, $2^2 - 1$, $2^3 - 1$, $2^4 - 1$ correct entries. This ``almost doubling the number of correct digits with each iteration" is quite characteristic of Newton's method. This clearly demonstrates quadratic convergence.

We then look at the secant method and Halley's method, which in continued fractions notation generate a Fibonacci number of convergents (by the secant method) and $\mathcal{O}(3^n)$ convergents showing cubic convergents (by Halley's method). This is not part of any standard numerical analysis curriculum, by the way, that we know of.

We have a fine opportunity for a digression on Fibonacci numbers; another for a digression on Schr{\"{o}}der iteration~\cite{schroder1993infinitely}.

\begin{remark}
	We avoid discussion of ``convergence" and note that each iterate solves $f(x) - f(x_n) = 0$; we \textsl{interpret} residuals in context.
\end{remark}

\subsection{Backward Error}
The reader will have noticed that we used $\sfrac{17}{12}$ as the exact square root of $2\ \sfrac{1}{144}$; similarly $\sfrac{577}{408}$ is the exact square root of $2\ \sfrac{1}{408^2}$. The is an example of what is called ``Backward Error Analysis.'' Instead of an approximate answer to the reference problem, we have an exact answer to a nearby problem.

This is a very powerful old idea, used by von Neumann, by Turing, and others. It was popularized in numerical analysis by Wilkinson. The book \cite{corless2013graduate} uses it in every chapter.

It is a useful idea pedagogically, as well due to its ease of understanding for most students. Students do understand it. If the course you're teaching the idea in has a formative assesment, the sly question comes up ``Is it all right if I give an answer to a slightly different question on the midterm?'' They get the idea.

Of course some problems are sensitive to changes in their formulation. This is also not quite mathematics: mathematics gains leverage by abstracting away context, and backward error needs problem context to be directly useful, but it can be done. If forward error (e.g.~$\left|x_n - \sqrt{2}\right|$) is needed, then one can introduce the notion of derivative and relative derivative (also known as logarithmic derivative or ``condition number''---von Neumann and Goldstine called it a ``figure of merit''~\cite{von1947numerical}).

But in this first year course, a direct interpretation allows one to avoid infinity and to avoid the calculus.
\begin{figure}[t]
	\centering
	\begin{tikzpicture}
		\draw (0, 0) -- (4, 0) -- (4, 4) -- (0, 0);
		\draw (2, 0) node[below]{$\frac{p_5}{q_5}$};
		\draw (4, 2) node[right]{$\frac{p_5}{q_5} = \frac{886731088897}{627013566048}$};
		\draw (1.95, 2) node[above, rotate=45]{2 + $\varepsilon$};
		\draw (3.8, 0) -- (3.8, 0.2) -- (4, 0.2);
	\end{tikzpicture}
	\caption{Example for backward error: instead of a 2 light year hypotenuse, we have a $2 + \varepsilon$ one, where $\varepsilon < 10^{-24}$ light years or $10^{-6}$ mm.}
	\label{fig:backward_error}
\end{figure}
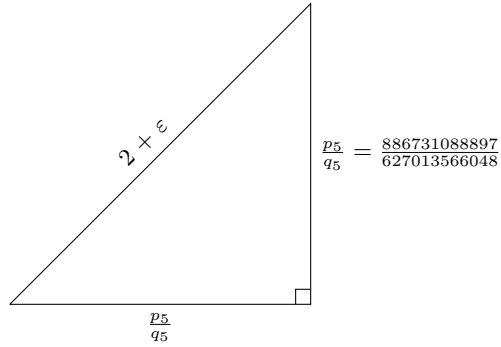
If $x = \sfrac{p_5}{q_5}$, then the difference between $x^2$ and $2$ is less than $2\cdot 10^{-24}$ light years or $\sim 10^{-6}$mm; that is, we've found the exact length of triangle whose hypotenuse is not $2$ light years, but 2 light years $+$ less than a millionth of a millimeter. Students get this.

The next section, on modified equations applies this idea in an advanced context for graduate students.

\subsection{Modified Equations}
 As an example lecture for the senior students, consider the construction of a so-called \textsl{modified equation}. This was done in lecture format. In an attempt to understand a very simple numerical method (fixed time-step Euler's method) on a very simple but nonlinear initial-value problem, namely $\dot y = y^2$, $y(0) = y_0$, one might attempt to use a technique called ``the method of modified equations". This leads to the series
\begin{equation}
	1 - v + \frac{3}{2}v^2 - \frac{8}{3}v^3 + \cdots - \frac{9427}{210}v^7 + \cdots \>.
\end{equation}
Computation of the series coefficients is somewhat involved but automatable. Consulting the Online Encyclopedia of Integer Sequences leads us to the work of Labelle~\cite{labelle1980inversion}, in combinatorics; this connection is powerful and unexpected.

The idea of this lecture is that we begin with ``brute force"; then consult the OEIS or some other resource to try to identify our results and find faster/better ways, and to make connections to other works. Then we extract other useful materials from the results, proving what we can. We give some more details of the lecture, below.

The method of modified equations is a neat idea: if we solve a DE, say $\dot y = f(y)$, by a numerical method, say $y_{n+1} = y_n + h \dot y_n$ (Euler's method), then in order to \textit{explain} the numerics we look for a modified equation
\begin{equation}
	\dot Y = F(Y;h)
\end{equation}
which (more nearly) has $Y(t_n) \doteq y_n$. References include \cite[pp.~662--671]{corless2013graduate}, and \cite{griffiths1986scope} as well as \cite{corless1994error}, and \cite{corless2016variations}. It's best explained by example. Consider $f(y) = y^2$. Then $\dot y = y^2$, $y(t_0) = y_0$ is actually easy to solve: $\frac{\dot y}{y^2} = 1$ so $\int^t_{t_0}\frac{\dot y \text{d}t}{y} = \int^t_{t_0} \text{d}t = t - t_0$ or $-\frac{1}{y}|^t_{t_0} = \frac{1}{y(t_0)} - \frac{1}{y(t)} = t - t_0$ or
\begin{equation}
	y(t) = \frac{y_0}{1 - y_0\cdot(t-t_0)} \>.
\end{equation}
This is singular at $t = t_0 + \frac{1}{y_0}$.

Let's pretend we don't know this, and instead solve the problem using Euler's method with a small fixed time step $n$:
\begin{equation}
	y_{n+1} = y_n + hy_n^2 \quad n = 0, 1, 2, \ldots\>.
\end{equation}
The key step in the method of modified equations is to replace this discrete recurrence relation with a functional equation: $Y(t)$ interpolating the data $(t_n, y_n)$ so $Y(t_n) = y_n$; since $t_{n+1} = t_n + h$ we have necessarily $Y(t_n+h) = Y(t_n) + hY(t_n)^2$.

We now insist that
\begin{equation}
	Y(t+h) = Y(t) + hY^2(t)
\end{equation}
everywhere, not just at $t=t_n$. This is a weird thing to do. Even weirder, we now look for a differential equation for $Y(t)$, to replace the (harder) functional equation. Using Taylor expansion we have
\begin{equation}
	Y(t+h) = Y(t) + \dot Y(t)h + \ddot Y(t)\frac{h^2}{2} + \dddot Y(t)\frac{h^3}{6} + \cdots
\end{equation}
so we have (in a statement with a breathtaking lack of rigour)
\begin{equation}
	Y(t) + h\dot Y(t) + \frac{h^2}{2} \ddot Y(t) + \frac{h^3}{6} \dddot Y(t) + \cdots = Y(t) + hY^2(t)
\end{equation}
or
\begin{equation}
	\dot Y(t) + \frac{h}{2} \ddot Y(t) + \frac{h^2}{3!} \dddot Y(t) + \cdots = Y^2(t) \>.
\end{equation}
This is a singular perturbation of the original. It's not very helpful. We can truncate:
\begin{equation}
	\sum^N_{k=1} Y^{(k)}(t)\frac{h^{k-1}}{k!} = Y^2(t) + \mathcal{O}(h^{N}) \>.
\end{equation}
We can differentiate (no breath left!)
\begin{equation}
	\ddot Y + \frac{h}{2} \dddot Y + \frac{h^2}{3!}Y^{\text{IV}} + \cdots = 2Y\dot Y
\end{equation}
and again
\begin{equation}
	\dddot Y + \frac{h}{2} Y^{\text{IV}} + \mathcal{O}(h^2) = 2(\dot Y)^2 + 2Y\ddot Y
\end{equation}
as many times as we need to.

Noticing that we need one less order in $h$ with each derivative we have
\begin{alignat}{4}
	\dot Y &{}+ \frac{h}{2} \ddot Y &{}+ \frac{h^2}{6} \dddot Y &= Y^2 + \mathcal{O}(h^3) \\
	\ddot Y &{}+ \frac{h}{2} \dddot Y & &= 2Y\dot Y + \mathcal{O}(h^2) \\
	\dddot Y & & &= 2(\dot Y)^2 + 2 Y \ddot Y + \mathcal{O}(h) \>.
\end{alignat}
Since the second equation implies $\ddot Y = 2 Y \dot Y + \mathcal{O}(h)$, and the first equation implies $\dot Y = Y^2 + \mathcal{O}(h)$, the third equation can be simplified to
\begin{align}
	\dddot Y &= 2Y^2 \dot Y + 2 Y (2 Y \dot Y) + \mathcal{O}(h) \nonumber \\
	&= 6 Y^2 \dot Y + \mathcal{O}(h).
\end{align}
Using this in the second equation gives
\begin{equation}
	\ddot Y + \frac{h}{2}(6 Y^2 \dot Y) = 2 Y \dot Y + \mathcal{O}(h^2)
\end{equation}
and using both of these in the first equation gives
\begin{equation}
	\dot Y + \frac{h}{2}(2Y \dot Y - 3hY^2 \dot Y) + \frac{h^2}{6}(6Y^2 \dot Y) = Y^2 + \mathcal{O}(h^3)
\end{equation}
or
\begin{equation}
	\dot Y + hY\dot Y - \frac{3}{2}h^2Y^2\dot Y + h^2Y^2 \dot Y = Y^2 + \mathcal{O}(h^3)
\end{equation}
or
\begin{equation}
	(1 + hY - \frac{1}{2}h^2Y^2)\dot Y = Y^2 + \mathcal{O}(h^3)
\end{equation}
or
\begin{equation}
	\dot Y = (1 + hY - \frac{1}{2}h^2Y^2)^{-1}Y^2 + \mathcal{O}(h^3)
\end{equation}
or
\begin{equation}
	\dot Y = (1 - hY + \frac{3}{2}h^2Y^2)Y^2 + \mathcal{O}(h^3) \>.
\end{equation}

What does this mean? It means that Euler's method applied to $\dot y = y^2$ gives a better solution to
\begin{equation}
	\dot y = (1 - hy + \frac{3}{2}h^2y^2)y^2
\end{equation}
(which wasn't intended); indeed it's $\mathcal{O}(h^3)$ accurate on \textit{that} one.

Why is this our first senior lecture on Open Problems for Experimental Maths? Because it was among RMC's first forays into the subject. He wrote a Maple program to compute more terms; then tried to use ``gfun" in Maple (which failed) and the Online Encyclopedia of Integer Sequences (which worked) to identify $B(v)$ where
\begin{equation}
	\dot y = B(hy)y^2
\end{equation}
with $B(v) = 1 - v + \frac{3}{2}v^2 + \cdots$ is the modification. From the OEIS we are led to \cite{labelle1980inversion} which gives (in another context)
\begin{equation}
	B(v) = \sum_{n\geq 0}c_nv^n
\end{equation}
with $c_0 = 1$ and $c_n = - \frac{1}{n}\sum^n_{i=1}\binom{n-i+2}{i+1}c_{n-i}$. This series \textit{diverges}. 
\begin{displayquote}
``Therefore, we may \textit{do} something with it" --- O. Heaviside~\cite{hardy1949divergent}.
\end{displayquote}
For very small $v$,
\begin{equation}
	B(v) \sim 1 - v + \frac{3}{2}v^2 - \frac{8}{3}v^3 + \mathcal{O}(v^4) \>.
\end{equation}	
Note: $y_{n+1} = y_n + hy_n^2 \ \rightarrow \ hy_{n+1} = hy_n + (hy_n)^2$ so putting $v(t) = hy(t)$ gives (with $\tau = \frac{t-t_0}{h}$)
\begin{align}
	v(\tau +1) &= v(\tau) + v^2(\tau) \\
	\frac{\text{d}v}{\text{d}\tau} &= B(v)v^2 \\
	\therefore \quad \frac{\text{d}v(\tau + 1)}{\text{d}\tau} &= B(v(\tau +1)\cdot v(\tau + 1)^2 \\
	\text{and} \quad \frac{\text{d}v(\tau + 1)}{\text{d}\tau} &= \frac{\text{d}v}{\text{d}\tau} + 2v\frac{\text{d}v}{\text{d}\tau} = (1 + 2v)\cdot B(v) \cdot v^2 \\
	\therefore \quad B(v+v^2)\cdot (v+v^2)^2 &= (1 + 2v)B(v)v^2 \\
	\text{or} \quad B(v+v^2)\cdot(1+v)^2 &= (1 + 2v)B(v) \quad \text{when} \ v \equiv 0 \>.
\end{align}
Therefore, we have a functional equation for $B$.
\begin{equation}
	B(v) = \frac{(1 + v)^2}{(1 + 2v)}\cdot B(v+v^2)
\end{equation}
or
\begin{align}
	B(v+v^2) &= \frac{(1 + 2v)}{(1 + v)^2} B(v) \\
	\frac{\text{d}v}{\text{d}t} &= B(v)v^2 \\
	\frac{\text{d}v/\text{d}t}{B(v)v^2} &= 1 \\
	\int^t_{t_0} \frac{\text{d}v/\text{d}\tau}{B(v)v^2}\text{d}\tau &= t - t_0
\end{align}
if we can evaluate $B(v)$, we can numerically integrate this to get $F(v(t)) - F(v(\tau_0)) = \tau - \tau_0$ which will allow us to plot $v(\tau)$.

What happens in the iteration $v_{n+1} = v_n + v_n^2$?
\begin{itemize}
	\item $v_n \rightarrow \infty$ as $ n \to \infty$ 
	\item $v_n \to 0$ as $n \to \infty$ 
	\item sometimes neither
\end{itemize}
We can see this \textit{experimentally}. Notice that $B$ has a \textit{pole} if $v = -\frac{1}{2}$ or if $v + v^2 = -\frac{1}{2}$ and so on. Therefore, the pre-images of $v_n = -\frac{1}{2}$ are \textit{all} poles. Therefore, poles approach 0 arbitrarily closely. Therefore, the series at 0 cannot converge.
\begin{figure}[h]
\centering
\includegraphics[width = 0.65\textwidth]{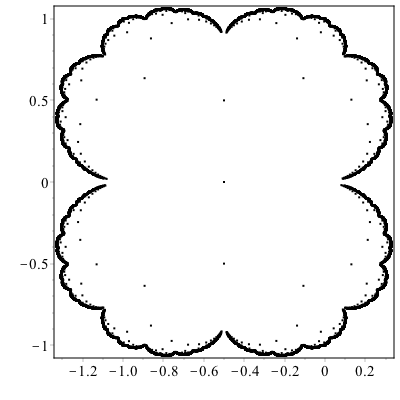}
\end{figure}

The Maple code to plot the pre-images is as follows:\\
\MapleInput{Digits := 15;}
\MapleOutput{Digits := 15}
\MapleInput{Wanted := 100000;}
\MapleOutput{Wanted := 100000}
\MapleInput{preimages := Array(0 ..Wanted);}
\MapleOutput{preimages := $\begin{bmatrix}\text{0 .. 100000 Array}\\ \text{Data Type: anything}\\ \text{Storage: rectangular}\\ \text{Order: Fortran\_order}\end{bmatrix}$}
\MapleInput{iworking := 0;}
\MapleOutput{iworking := 0}
\MapleInput{ifree := 1;}
\MapleOutput{ifree := 1}
\MapleInput{preimages[0] := 0.5;}
\MapleOutput{preimages\textsubscript{0} := -0.5}
\MapleInput{while ifree $<$ Wanted do\\
\hspace*{1cm} p := preimages[iworking];\\
\hspace*{1cm} preimages[ifree] := $-\dfrac{(1 + \sqrt{1 + 4\cdot p})}{2}$; \\
\hspace*{1cm} ifree := ifree + 1;\\
\hspace*{1cm} preimages[ifree] := $-\dfrac{p}{\text{preimages[ifree - 1]}}$;\\
\hspace*{1cm} ifree := ifree + 1;\\
\hspace*{1cm} iworking := iworking + 1;\\
end do:}\\
\MapleInput{plot(map(t$\to$[Re(t), Im(t)], [seq(preimages[k], k = 0 ..Wanted)]), \\
style=point, colour=black, symbol=point, symbolsize=1)}\\

Remarks
\begin{itemize}
	\item The pre-images of $-\frac{1}{2}$ are a subset of the Julia set. $B(v)$ is singular on \textit{all} of those!
	\begin{equation}
		B(p_N) = \prod^{N-1}_{j=0} \frac{(1 + v_j)^2}{(1 + 2v_j)}B(v_N) = B(-\frac{1}{2}) = \infty
	\end{equation}	
	\item The set of pre-images of $-\frac{1}{2}$ is an infinite set and is a \textsl{fractal} (looks like a cauliflower) and approaches 0 arbitrarily closely. Therefore, there are poles of $B(v)$ arbitrarily close to 0 and thus, the series $B(v) = 1 - v + \cdots$ \textsl{cannot} converge.
	\item The pre-images of $-1$ are \textsl{zeros} of $B$ and they're also dense in the Julia set. Therefore, the Julia set forms a \textsl{natural boundary} for $B(v)$ and $\frac{\text{d}v}{\text{d}t} = B(v)v^2$ is \textsl{nasty.}
	\item The pre-images for the Mandelbrot derivation are computed similarly but more simply.
	\item We re-used the pre-image code for the first-year course, for general Julia sets.
\end{itemize}

\subsection{Chaos Game Representation of DNA sequences}
This final vignette shows that these dynamical idea are not just toys or idle puzzles but useful tools for science. As mentioned in the previous section, when the course was first taught, this research area was new and not yet published. Since then, \cite{karamichalis2016modmaps3d}, \cite{karamichalis2015investigation}, and \cite{kari2015mapping} have been published, in which the most recent paper (June 2017) was published in a top journal, Bioinformatics. 

It is general knowledge that DNA is a double helix and is made of four bases: adenine (A), guanine (G), cytosine (C), and thymine (T). However, for DNA sequencing, we only need to look at one of the two strands of the DNA since the bases are paired: A with G, and C with T. This means that if we know one strand, we automatically would know both strands (as long as there are no mutations). DNA sequences are unique to each organism; therefore, using chaos game representation, the goal is to be able to quantitatively classify the organisms by their classes (in this case, animal classes; the five most well-known classes of vertebrates are mammals, birds, fish, reptiles, and amphibians).

\begin{figure}[h]
	\centering
	\begin{tikzpicture}
		\filldraw (-2, -2) circle (2pt) node[below]{(-1, -1)} node[above]{A};
		\filldraw (-2, 2) circle (2pt) node[below]{(-1, 1)} node[above]{C};
		\filldraw (2, -2) circle (2pt) node[below]{(1, -1)} node[above]{T};
		\filldraw (2, 2) circle (2pt) node[below]{(1, 1)} node[above]{G};
	\end{tikzpicture}
	\caption{Chaos game representation setup}
	\label{fig:chaos}
\end{figure}
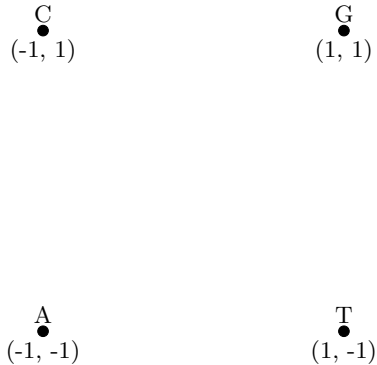

We start with a plot with four corners, each of which represents one of the four bases in DNA, as shown in figure \ref{fig:chaos}. The algorithm for chaos game representation is outlined in the following steps.
\begin{enumerate}
	\item Put a dot in the center; this is the ``current dot" at the beginning.
	\item Pick the next letter in the sequence (if none, stop), and draw a invisible line from the current dot towards the corner representing the base.
	\item Put a dot halfway on the line; this becomes the new ``current dot."
	\item Return to step 2
\end{enumerate}
The short Maple code below explains this precisely.

However, for the course, we did not use any DNA sequences; instead, we generated a sequence of (not so) random numbers, in which each number represents an index of our variable \texttt{ACGT}, which corresponds to a coordinate point (representing a base). The following is the Maple  procedure that was used to generate chaos game representation plots.\\
\MapleInput{CGR := proc(s::list)\\
\hspace*{1cm}local ACGT, k, n, p;\\
\hspace*{1cm}ACGT := table(); \\
\hspace*{1cm}ACGT[0] := [-1, -1]; \\
\hspace*{1cm}ACGT[1] := [-1, 1]; \\
\hspace*{1cm}ACGT[2] := [1, 1]; \\
\hspace*{1cm}ACGT[3] := [1, -1]; \\
\hspace*{1cm}n := nops(s);\\
\hspace*{1cm}p := Array(0 ..n);\\
\hspace*{1cm}p[0] := [0, 0]; \\
\hspace*{1cm}for k to n do\\
\hspace*{1.75cm}p[k] := (ACGT[s[k]] + p[k-1])/2.0; \\
\hspace*{1cm}end do: \\
\hspace*{1cm}plots[pointplot]([seq(p[k], k=1..n)], colour=black, axes=none); \\
end proc:}\\
We can see that in the code, instead of using letters for the bases, we have used the numbers 0, 1, 2, 3 instead.

\begin{figure}
	\centering
	\begin{subfigure}[t]{0.45\textwidth}
		\includegraphics[width=\textwidth]{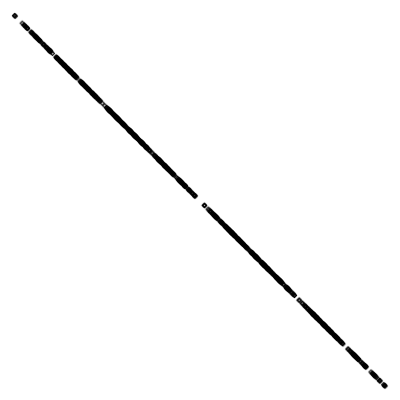}
		\caption{$103 + k$th prime mod 4}
		\label{fig:cgr_ithprime}
	\end{subfigure}
	\qquad
	\begin{subfigure}[t]{0.45\textwidth}
		\includegraphics[width=\textwidth]{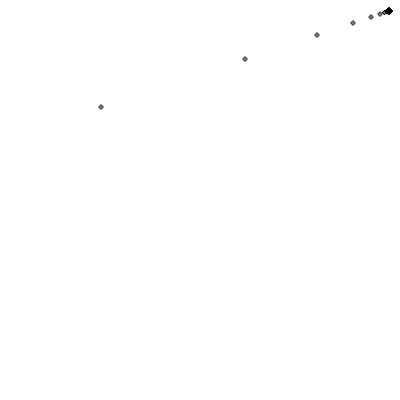}
		\caption{Partial quotients of continued fraction of $\sqrt{2}$}
		\label{fig:cgr_cf_sqrt_2}
	\end{subfigure}

	\begin{subfigure}[t]{0.45\textwidth}
		\includegraphics[width=\textwidth]{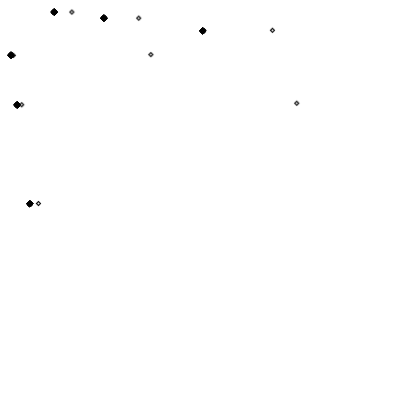}
		\caption{Partial quotients of continued fraction of $e$}
		\label{fig:cgr_exp}
	\end{subfigure}
	\qquad
	\begin{subfigure}[t]{0.45\textwidth}
		\includegraphics[width=\textwidth]{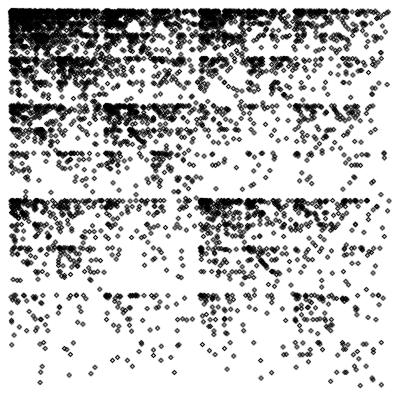}
		\caption{Partial quotients of continued fraction of $\pi$}
		\label{fig:cgr_cf_pi}
	\end{subfigure}
	
	\begin{subfigure}[t]{0.45\textwidth}
		\includegraphics[width=\textwidth]{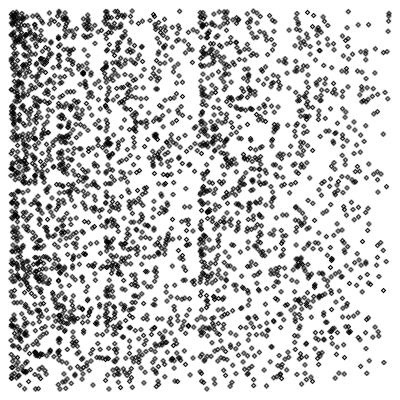}
		\caption{Digits of $\pi$}
		\label{fig:cgr_digits_pi}
	\end{subfigure}
	\qquad
	\begin{subfigure}[t]{0.45\textwidth}
		\includegraphics[width=\textwidth]{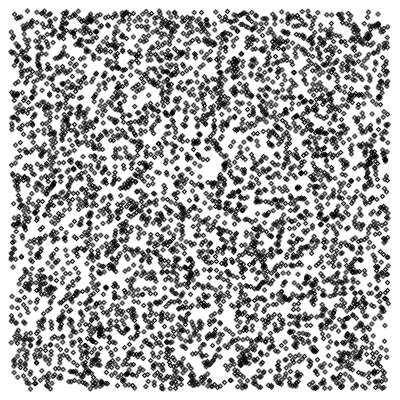}
		\caption{Random numbers}
		\label{fig:cgr_random}
	\end{subfigure}
	
	\caption{Examples of chaos game representations (all examples are mod 4)}
	\label{fig:cgr}
\end{figure}

The first sequence the class looked at was the $i$th prime number (actually, we took the $103+k$th prime). The result, a diagonal line (figure \ref{fig:cgr_ithprime}), was surprising at first. But in actually, it is not \textsl{that} surprising in hindsight. As we all know, prime numbers bigger than 2 will never be even, so the coordinate points (in this case $(1, 1)$ and $(-1, -1)$) that represented ``even'' values would never occur, thus creating a straight diagonal line. If we had reassigned the values the coordinate points represent, the plot would turn out different; instead of a diagonal line, it could possibly be a horizontal or vertical line. One needs to be mindful of how the coordinate points are assigned.

The students were then encouraged to try other sequences. Some started with something familiar: one recurring theme of this course was $\sqrt{2}$, so with this in mind, some students jumped on the opportunity to plot the chaos game representation of the quotients of the continued fraction of $\sqrt{2}$. As we have seen earlier in the course, the sequence goes like
\begin{equation}
	1 + [2, 2, 2, 2, \ldots, 2] \>.
\end{equation}
Because all elements of the sequence (apart from the first one) are 2's, it is not surprising to see a (faint) diagonal line with most of the points in the upper right corner, shown in figure \ref{fig:cgr_cf_sqrt_2}. The students also tried the partial quotient of $e$, shown in figure \ref{fig:cgr_exp}. What is seen here makes sense as the sequence mostly contains 1's and these alternates with even values, so in this case, either 0 or 2. Therefore, it is clear that there are no points in the lower right corner since that represents the value 3. Unfortunately, these two plots do not look all that impressive; in fact, it is pretty underwhelming. Disappointed with this result, the students decided to experiment with other sequences in which all four coordinates occur.

The students then thought, ``Why not take the partial quotients of the continued fraction of $\pi$? The results \textsl{must} be random.'' As shown in figure \ref{fig:cgr_cf_pi}, unexpectedly, there is a indeed pattern, which shows that the sequence of partial quotients of the continued fraction of $\pi$ is not as random as we thought at the beginning of the course. Although the continued fraction quotients of $\pi$ may not actually be random, it is known that the digits of $\pi$ themselves are. The chaos game representation plot (figure \ref{fig:cgr_digits_pi}) of this shows some randomness. This can be compared to the chaos game representation plot of randomly-generated values mod 4 (figure \ref{fig:cgr_random}).

Additionally, this brings up the thought of whether this idea can be extended to other figures---base 6, base 8, base 3. This could possibly be explored in a future version of the course. This connects interestingly with the random walk on the digits of $\pi$~\cite{artacho2013walking}.

\section{Concluding Remarks}
There were many successful outcomes from this course. The most political achievement was meeting the Hon.~Kathleen Wynne, Premier of Ontario, at the time of publication as part of a WALS demo (she said that she was frightened of mathematics). Secondly, two undergraduate students who participated in the first year of the course (Alex Wu and Tiam Koukpari) founded Mustang Capital, initially known as The Algorithm Trading Club, which has over 100 members. They use the pedagogical principles of this course to teach themselves about algorithmic trading. There were also many academic achievements that came from this course: beautiful and prize-winning posters (in ISSAC 2016) were created, and two students (Yang Wang and Ao Li) published papers from their projects. In addition to this, RMC recruited at least one PhD student (possibly more are pending). Lastly, interest in the course has been generated for The Digital Humanities Program.

\subsection*{A less happy outcome}
The course AM 1999 was taught only once, in spite of the promises from the outgoing Chair and the Dean to run it for two years. RMC had failed to get the Dean's promise in writing, and the incoming acting Chair cancelled the second offering, owing to an important misunderstanding on the part of some colleagues. Jon Borwein was scathing about this decision---any new program needs time for growth of awareness, and this was especially acute here given the short period from conception (May 2014) to first delivery (Sept 2014). The senior course AM 9619 \textsl{was} offered a second time, which would have made the second offering of AM 1999 ``free"; this makes the lost opportunity even more sharp.

This bit of data is not included here as a lament, but rather for clarity and as a recommendation: before undertaking such a serious undertaking, get your promised support in writing, and be sure to tell your colleagues what you are doing.

\section*{Acknowledgements}
We thank David Jeffrey for his early encouragement, Steven Thornton and Torin Viger for their help with the course material, and the membership of USAM: Julia Jankowski, Andy Wilmot, and Anna Zhu. Gavin Watson, Stephanie Oliver, and Wendy Crocker were very helpful with active learning and the use of WALS. We thank Lila Kari for introducing us to the chaos game representation. We also thank the Rotman Institute of Philosophy and the Fields Institute for Research in the Mathematical Sciences for their sponsorship of the Computational Discovery conference (\url{acmes.org}), at which some of these ideas presented here were refined.

\bibliographystyle{spmpsci}
\bibliography{Experimental-math}

\end{document}